\documentclass[12pt]{article}
\usepackage{latexsym}
\usepackage{graphicx}
\usepackage[dvips]{epsfig}
\usepackage{amsfonts}
\usepackage{amssymb}
\newtheorem{biggy}{Theorem}

\newtheorem{riggy}{Remark}
\newtheorem{liggy}{Lemma}
\newtheorem{piggy}{Proposition}
\newcommand{\hole}{{\mkern2mu|\mkern2mu}}
\begin{document}

\title{Petrie-Coxeter Maps Revisited}

\author{Isabel Hubard, Egon Schulte\thanks{Supported by NSA-grant H98230-05-1-0027} 
\ and Asia Ivi\'{c} Weiss\thanks{Supported by NSERC of Canada Grant \#8857}}

\maketitle


\begin{abstract}
\noindent
This paper presents a technique for constructing new chiral or regular polyhedra (or
maps) from self-dual abstract chiral polytopes of rank $4$. From improperly self-dual chiral
polytopes we derive ``Petrie-Coxeter-type" polyhedra (abstract chiral analogues of the
classical Petrie-Coxeter polyhedra) and investigate their groups of automorphisms.
\end{abstract}

\section{Introduction}

The modern theory of regular polytopes and their geometric realizations in Euclidean spaces has 
been greatly influenced by the discovery of the Petrie-Coxeter polyhedra in $\mathbb{E}^3$
and Coxeter's regular skew polyhedra in $\mathbb{E}^4$ (see Coxeter~\cite{skew},
Dress~\cite{dress}, Gr\"unbaum~\cite{gruen}, McMullen~\cite{full,fourdim}, and
\cite{AMSpaper,arp}).  Two of the Petrie-Coxeter polyhedra and all of Coxeter's
skew polyhedra inherit their symmetry properties directly from self-dual regular figures,
namely the cubical tessellation in $\mathbb{E}^3$, and the self-dual regular convex
polytopes (the $4$-simplex $\{3,3,3\}$ and $24$-cell $\{3,4,3\}$) or the double $p$-gonal
prisms (or $\{p,2,p\}$) in $\mathbb{E}^4$, respectively (see Section~\ref{chimprop}).
However, the symmetries of these figures provide only one half of the symmetries of the
polyhedra; the other half corresponds to the dualities (see \cite{peter-egon,
schulte-wills}). 

Coxeter credited Petrie with the original idea of admitting infinite regular polyhedra in 
$\mathbb{E}^3$ with skew vertex-figures. In particular, Petrie discovered two of the
Petrie-Coxeter polyhedra, and then Coxeter found the third and carried out the corresponding
enumeration in $\mathbb{E}^4$. In \cite{garner}  and \cite{weiss-lucic}, this enumeration is
extended to hyperbolic $3$-space $\mathbb{H}^3$, where there are thirty-one regular skew
polyhedra.

The connection between self-dual polytopes and skew polyhedra was further investigated in
\cite{peter-egon} in the context of abstract regular polytopes, where a construction of
regular ``Petrie-Coxeter-type" polyhedra from self-dual regular $4$-polytopes is described.
In this paper we extend these ideas to abstract chiral polytopes. We exploit the two
different kinds of self-duality (proper and improper self-duality) for chiral polytopes and
present constructions of polyhedra, which are chiral or regular, depending on the kind of
self-duality of the original chiral $4$-polytope. From improperly self-dual chiral
$4$-polytopes we obtain chiral ``Petrie-Coxeter-type" polyhedra (or maps), that is,
abstract chiral analogues of the ordinary Petrie-Coxeter polyhedra.

\section{Regular and chiral polytopes}

In this section we briefly review some definitions and basic results from the theory of
abstract polytopes. For details we refer to \cite{arp}.

A {\em polytope of rank\/} $n$, or an {\em $n$-polytope\/}, is a partially ordered set $\cal P$ 
with a strictly monotone rank function having range $\{-1,0, \ldots, n\}$. The elements of $\cal
P$ with rank $j$ are the {\em $j$-faces\/} of $\cal P$. The faces of ranks $0$, $1$ or $n-1$
are called {\em vertices\/}, {\em edges\/} or {\em facets\/}, respectively. There is a smallest
face $F_{-1}$ (of rank $-1$) and a greatest face $F_n$ (of rank $n$), and each {\em flag\/}
(maximal chain) contains exactly $n+2$ faces. Also,  $\cal P$ is {\em strongly
flag-connected\/}, that is, any two flags $\Phi$ and $\Psi$ of $\cal P$ can be joined by a
sequence of flags $\Phi = \Phi_0, \Phi_1, \ldots, \Phi_k = \Psi$ such that each $\Phi_{i-1}$
and $\Phi_i$ are {\em adjacent\/} (they differ by just one face), and $\Phi \cap \Psi \subseteq \Phi_i$ for each $i$. Furthermore, $\cal P$ has the following homogeneity property:\  whenever  $F \leq G$, with $\textrm{rank}(F) = j-1$ and $\textrm{rank}(G) = j+1$, then there are exactly two $j$-faces H with $F \leq H \leq G$. These conditions essentially say that $\cal P$ shares many combinatorial properties with classical polytopes.

If $\Phi$  is a flag, then by $\Phi^i$ we denote the unique flag adjacent to $\Phi$ and differing 
from $\Phi$ in the face of rank $i$. Note that $(\Phi^i)^i = \Phi$. We extend our notation
and write $\Phi^{i_{1},\ldots,i_{k-1},i_{k}} := (\Phi^{i_{1},\ldots,i_{k-1}})^{i_{k}}$.  For
any two faces $F$ and $G$ with $F \leq G$, we call $G/F :=\{ H \mid F \leq H \leq G \}$ a {\em section\/} of $\cal P$; if $F$ is a vertex and $G=F_{n}$, then this is the {\em
vertex-figure\/} of $\cal P$ at $F$.  Note that $G/F$ is an abstract polytope of rank equal to
$\textrm{rank}(G) - \textrm{rank}(F) - 1$.

The group of automorphisms (order preserving bijections) of $\cal P$ is denoted by 
$\Gamma({\cal P})$. If $\cal P$ admits an order reversing bijection, then $\cal P$ is said to
be {\em self-dual\/} and any such bijection is called a {\em duality\/} of $\cal P$. In this
case the group of all automorphisms and dualities is called the {\em extended group\/} of $\cal P$ and is denoted by $\overline{\Gamma}({\cal P})$.  A {\em polarity\/} is an involutory duality.  Note that, if $\Phi$ is any flag and $\delta$ any duality, then $\Phi^j\delta = (\Phi\delta)^{n-j-1}$ for $j=0,\ldots,n-1$.

A polytope $\cal P$ is {\em regular\/} if $\Gamma({\cal P})$ is transitive on the flags of $\cal P$.
The group of a regular polytope $\cal P$ is generated by involutions $\rho_i$, $i=0,\ldots,n-1$, mapping a fixed, or {\em base\/}, flag $\Phi$ to the adjacent flags $\Phi^i$.
These  generators satisfy (at least) the relations 
\begin{equation}
\label{equatone}
(\rho_i \rho_j)^{p_{ij}}=\epsilon \;\; \textrm{ for } i,j=0, \ldots,n-1, 
\end{equation}
where $p_{ii}=1$, $p_{ji} = p_{ij} =: p_{i+1}$ if $j=i+1$, and $p_{ij}=2$ otherwise. We then say 
that ${\cal P}$ is of {\em type\/} $\{ p_1,\ldots,p_{n-1} \}$. Furthermore, these generators
satisfy the {\em intersection condition\/},
\begin{equation}
\label{equattwo}
\langle \rho_i \mid i \in I \rangle \cap \langle \rho_i \mid i \in J \rangle 
= \langle \rho_i \mid i \in {I \cap J} \rangle 
\quad (I,J \subseteq \{0,1,\ldots,n-1\}). 
\end{equation}

Conversely, if $\Gamma = \langle \rho_0, \ldots, \rho_{n-1} \rangle$ is a group whose generators 
satisfy relations (\ref{equatone}) and condition (\ref{equattwo})\  (such groups are called
{\em C-groups}), then there exists a regular polytope with $\Gamma({\cal P}) \cong \Gamma$.

We note that every self-dual regular polytope possesses a polarity mapping the base flag to itself (but reversing the order of its elements). 

A polytope $\cal P$ (of rank $n \geq 3$) is {\em chiral\/} if $\cal P$ is not regular, but if for 
some base flag $\Phi = \{F_{-1}, F_0, \ldots, F_n \}$ there exist automorphisms $\sigma_1,
\ldots, \sigma_{n-1}$ of $\cal P$ such that $\sigma_i$ fixes all faces in $\Phi \setminus \{
F_{i-1}, F_{i} \}$ and cyclically permutes consecutive $i$-faces of $\cal P$ in the section
$F_{i+1} / F_{i-2} = \{G \mid F_{i-2} \leq G \leq F_{i+1} \}$ of rank $2$. We may choose
$\sigma_i$ in such a way that 
\[ \Phi\sigma_i = \Phi^{i,i-1},\; \textrm{ for all } i \in \{1,\ldots,n-1\} ,\]
and hence
\[ \Phi\sigma_i^{-1} = \Phi^{i-1,i}, \textrm{ for all } i \in \{1,\ldots,n-1\} . \]
The automorphisms $\sigma_1,\ldots,\sigma_{n-1}$ generate $\Gamma({\cal P})$ and satisfy (at
least) the relations
\begin{equation}
\label{equatthree}
\begin{array}{rcl}
\sigma_i^{p_i}                          & = & \epsilon\;\, \textrm{ for } 1\leq i \leq n-1, \\[.03in]
(\sigma_i \sigma_{i+1}\dots \sigma_j)^2 & = & \epsilon\;\, \textrm{ for } 1\leq i < j \leq n-1,
\end{array} 
\end{equation}
where again $\{ p_1,\ldots,p_{n-1} \}$ is the {\em type} of $\cal P$. For general background on chiral polytopes see \cite{chiral,chirality}.

Alternatively, chiral polytopes can be characterized as those polytopes whose automorphism groups have two flag orbits, with adjacent flags in distinct orbits (see \cite {chiral, chirality}).
Consequently, if a chiral polytope $\cal P$ is self-dual, it can be so in one of two ways. If a
duality of $\cal P$ maps a flag to a flag (in reverse order) in the same orbit under $\Gamma({\cal P})$, then each duality must do so with each flag; in this case $\cal P$ is said to be {\em properly self-dual\/}. However, if a duality maps a flag to a flag in a different orbit under $\Gamma({\cal P})$, then again each duality must do so with each flag, and in this case $\cal P$ is said to be {\em improperly self-dual\/}.

For a regular $n$-polytope $\cal P$ with group $\Gamma({\cal P}) =
\langle \rho_{0},\rho_{1},\ldots,\rho_{n-1} \rangle$, define $\sigma_{i} :=
\rho_{i-1}\rho_{i}$ for $i=1,\ldots,n-1$. Then  $\sigma_{1},\ldots,\sigma_{n-1}$ generate the
({\em combinatorial\/}) {\em rotation subgroup\/} $\Gamma^{+}({\cal P})$ of $\Gamma({\cal P})$ (of
index at most $2$) and have properties similar to those of the distinguished generators for
chiral polytopes.

From now on we restrict our considerations to polytopes of rank $n \leq 4$. The generators
$\sigma_i$ of a chiral polytope $\cal P$ also satisfy an {\em intersection condition\/}, which for
$n=4$ is
\begin{equation}
\label{intchir4}
\begin{array}{c}
\langle \sigma_1, \sigma_2 \rangle \cap \langle  \sigma_2, \sigma_3 \rangle 
= \langle \sigma_2 \rangle,\\[.05in]
\langle\sigma_1\rangle \cap \langle\sigma_2\rangle = \{\epsilon\} = 
\langle \sigma_2 \rangle \cap \langle \sigma_3 \rangle,
\end{array}
\end{equation}
while for $n=3$ we only have 
\begin{equation}
\label{intchir3}
\langle \sigma_1  \rangle \cap \langle \sigma_2 \rangle = \{\epsilon\}.
\end{equation}
Conversely, if $\Gamma  = \langle \sigma_1, \ldots, \sigma_n \rangle$ is a group whose
generators satisfy the relations (\ref{equatthree}) and this intersection condition, namely
(\ref{intchir4}) for $n=4$ or (\ref{intchir3}) for $n=3$, then there exists a chiral polytope 
$\cal P$ with $\Gamma({\cal P}) \cong \Gamma$ or a regular polytope $\cal P$ with 
$\Gamma^{+}({\cal P}) \cong \Gamma$.

\section{Holes, zig-zags, and Petrie-polygons}

Each $3$-polytope $\cal P$ gives rise to a map on a surface, which in the case of a chiral
polytope must be orientable. A {\em $j$-hole\/} of $\cal P$ is an edge-path on the surface
which leaves a vertex by the $j$-th edge from which it entered, always in the same sense (in
some local orientation), keeping to the left (say). Similarly, a {\em $j$-zigzag\/} of $\cal P$
is an edge-path on the surface which leaves a vertex by the $j$-th edge from which it entered,
alternating the sense. For example, the $1$-holes are simply the $2$-faces of $\cal P$, and the
$1$-zigzags are precisely the {\em Petrie polygons\/} of $\cal P$ (see \cite{cm}). 

For a regular $3$-polytope $\cal P$, the length of a $j$-hole or $j$-zigzag is given by the
period of $\rho_0\rho_1(\rho_2\rho_1)^{j-1}$ or $\rho_0 (\rho_1 \rho_2)^j$ in $\Gamma({\cal
P})$, respectively (see \cite[pp. 196-197]{arp});  note here that
$\rho_0\rho_1(\rho_2\rho_1)^{j-1} = \sigma_1 \sigma_2^{1-j}\in \Gamma^{+}({\cal P})$
(but $\rho_0 (\rho_1 \rho_2)^j \in \Gamma^{+}({\cal P})$ only if $\Gamma^{+}({\cal P}) =
\Gamma({\cal P})$).  Similarly, for a chiral $3$-polytope $\cal P$, the length of a $j$-hole is the period of $\sigma_1 \sigma_2^{1-j}$ in $\Gamma({\cal P})$. 

If a regular $3$-polytope of type $\{p, q\}$ is completely determined by the lengths $h_j$ of its
$j$-holes, for $2 \leq j \leq  k:= \lfloor q/2 \rfloor$, and lengths $r_j$ of its $j$-zigzags,
for $1 \leq j \leq  k:= \lfloor q/2 \rfloor$, then we denote it by
\[ {\cal P} = \{p,q \hole h_2, \ldots, h_k \}_{r_1, \ldots,r_k},  \]
with the convention that any unnecessary $h_j$ or $r_j$ (that is, one not needed for the 
specification) is replaced by a $\cdot$, with those at the end of the sequence omitted. More
generally, we say that a regular $3$-polytope is of {\em type\/} 
$\{p,q \hole h_2, \ldots, h_k \}_{r_1,\ldots,r_k}$ if its $j$-holes and $j$-zigzags have
lengths $h_j$ or $r_j$, respectively, with $j$ as above; we also use this terminology for
chiral $3$-polytopes when applicable (that is, when only the $h_j$'s are specified).

For chiral or regular $4$-polytopes $\cal P$ we define two particular automorphisms 
\[ \pi_L := \sigma_1\sigma_3 \;\, \mbox{ and } \;\, \pi_R := \sigma_1\sigma_3^{-1} . \] 
Let $F_0$ and $F_1$ be the base vertex and base edge (of the base flag $\Phi$), respectively.  Then the orbits $F_0 \langle \pi_L \rangle$ and $F_1\langle \pi_L\rangle$, respectively,  are the vertex-set and edge-set of a polygon; its images under $\Gamma({\cal P})$ are called the
{\em left-handed Petrie polygons}, or for short {\em left-Petrie polygons\/}, of $\cal P$.  Similarly, we obtain {\em right-Petrie polygons} from $\pi_R$. Generally, for chiral polytopes, $\pi_L$ and
$\pi_R$ have different orders, so the left- and right-Petrie polygons have different
lengths. Left- and right-Petrie polygons were employed by Coxeter in \cite{th} to construct
chiral (or, in his terminology, {\em twisted\/}) $4$-polytopes. For regular $4$-polytopes, the elements $\pi_L$ and $\pi_R$ are conjugate in $\Gamma({\cal P})$, so left- and right-Petrie polygons necessarily have the same lengths. 

\section{Chiral maps from improperly self-dual chiral polytopes}
\label{chimprop}

From \cite{peter-egon} we know that with any self-dual regular $4$-polytope $\cal P$ of type
$\{p,q,p\}$ is associated a certain regular map of type $\{4,2q \hole p\}$ (with group
$\overline{\Gamma}({\cal P})$), called the {\em Petrie-Coxeter polyhedron of $\cal P$\/}; its
dual is of type $\{2q,4 \hole  p\}$.  When $\cal P$ is the universal polytope $\{p,q,p\}$ (see 
\cite{arp}), we obtain (the universal) $\{4,2q \hole p\}$ (as well as its dual $\{2q,4\hole 
p\}$). For example, from the $4$-simplex $\{3,3,3\}$ and $24$-cell $\{3,4,3\}$ in Euclidean
$4$-space $\mathbb{E}^4$ we derive Coxeter's finite skew polyhedra $\{4,6 \hole 3\}$ and
$\{4,8 \hole 3\}$ in $\mathbb{E}^4$ (see \cite{skew}, and for figures of projections into
$\mathbb{E}^3$ see \cite{schulte-wills}), while the regular tessellations $\{4,3,4\}$ and
$\{5,3,5\}$ of Euclidean or hyperbolic $3$-space $\mathbb{E}^3$ or $\mathbb{H}^3$,
respectively, yield the (classical) Petrie-Coxeter polyhedron $\{4,6 \hole 4\}$ in
$\mathbb{E}^3$ and the hyperbolic skew polyhedron $\{4,6 \hole 5\}$ in $\mathbb{H}^3$ (see
\cite{skew,garner,weiss-lucic}). The toroidal regular skew polyhedra $\{4,4\hole p\}$ in
$\mathbb{E}^4$ are related in a similar fashion to the double $p$-gonal prism
$\{p\}\times\{p\}$ (or $\{p,2,p\}$ in $\mathbb{S}^3$).

The construction of \cite{peter-egon} employs the self-duality of $\cal P$. Since $\cal P$ is
regular, it has a (unique) polarity $\omega$ that fixes the base flag and thus permutes (by conjugation) the generators of $\Gamma({\cal P}) = \langle \rho_{0},\ldots,\rho_{3} \rangle$ according to $\omega\rho_{i}\omega = \rho_{3-i}$, with $i=0,1,2,3$. For the universal polytope $\{p,q,p\}$, this polarity corresponds to the symmetry of the underlying (string) Coxeter diagram. For an arbitrary (self-dual and regular) polytope $\cal P$, the generators $\tau_{0},\tau_{1},\tau_{2}$ of the automorphism group $\overline{\Gamma}({\cal P})$ of its Petrie-Coxeter polyhedron are obtained from a twisting operation (in the sense of \cite[Ch.8]{arp}) on $\Gamma({\cal P})$, 
\begin{equation}
\label{twistpcp}
(\rho_0,\rho_1,\rho_2,\rho_3,\omega) \mapsto 
(\rho_0,\omega,\rho_2) =: (\tau_{0},\tau_{1},\tau_{2}) .
\end{equation}
\smallskip
Alternatively we may view (\ref{twistpcp}) as a mixing operation (in the sense of
\cite[Ch.7]{arp}) on $\overline{\Gamma}({\cal P})$.

We now make use of the fact that the automorphism groups of chiral polytopes of rank greater
than $3$ are generated by involutions, and apply similar techniques to construct new
``Petrie-Coxeter type" maps from self-dual {\em chiral\/} polytopes. In this Section we
discuss improper self-duality.

Henceforth, we assume that $\cal P$ is a chiral $4$-polytope of type $\{p,q,p\}$, $\Phi$ its
base flag, and $\sigma_{1},\sigma_{2},\sigma_{3}$ the distinguished generators of
$\Gamma({\cal P})$ (with respect to $\Phi$). Then these generators satisfy (at least) the
following relations
\begin{equation}
\sigma_1^p = \sigma_2^q = \sigma_3^p = 
(\sigma_1 \sigma_2)^2 = (\sigma_1 \sigma_2 \sigma_3)^2 = (\sigma_2 \sigma_3)^2 = \epsilon.
\end{equation}
Note that $\sigma_1\sigma_2, \sigma_1\sigma_2\sigma_3, \sigma_2 \sigma_3$ is a set of
involutory generators of $\Gamma({\cal P})$.

Let $\cal P$ be improperly self-dual. Then there exists a duality $\delta$ which exchanges
the two flag orbits. In fact, by multiplying $\delta$ with an automorphism of $\cal P$ (if need
be), we may choose $\delta$ so that $\Phi\delta = \Phi^{0}$. Then,
\[ \Phi^{0}\delta = \Phi^{0,3},\;\; \Phi^{0,3}\delta = \Phi^{3},\;\; 
\Phi^{3}\delta = \Phi\, . \]
Bear in mind here that $\delta$, being a duality, maps a $j$-face in a flag to the $(3-j)$-face in the image flag. It follows that the {\em distinguished duality\/} $\delta$ has period $4$ and satisfies the
relations 
\begin{equation}
\label{delacts}
\begin{array}{rcl}
\delta^{2}                  &=& \sigma_{1}\sigma_{2}\sigma_{3} , \\
\delta^{-1}\sigma_{1}\delta &=& \sigma_{3}^{-1} , \\
\delta^{-1}\sigma_{2}\delta &=& \sigma_{1}\sigma_{2}\sigma_{1}^{-1}, \\
\delta^{-1}\sigma_{3}\delta &=& \sigma_{1} , 
\end{array} 
\end{equation}
in $\overline{\Gamma}({\cal P})$ (see \cite[Theorem 3.5]{self-duality}). Moreover, it is easy
to see that $\delta$ fixes the basic involution $\sigma_{1}\sigma_{2}\sigma_{3}$, that is,
$\delta^{-1} \sigma_{1}\sigma_{2}\sigma_{3} \delta = \sigma_{1}\sigma_{2}\sigma_{3}$,
while 
\begin{equation}
\label{moredel}
\begin{array}{rcl}	
\delta^{-1}\sigma_{1}\sigma_{2}\delta &=& \sigma_{1}\sigma_{2}\sigma_{3}\sigma_{1}^{-1},\\
\delta^{-1}\sigma_{1}\sigma_{2}\sigma_{3}\sigma_{1}^{-1}\delta &=& 
\sigma_{3}^{-1}\sigma_{1}\sigma_{2}\sigma_{3},\\
\delta^{-1}\sigma_{3}^{-1}\sigma_{1}\sigma_{2}\sigma_{3}\delta &=& \sigma_{2}\sigma_{3},\\
\delta^{-1}\sigma_{2}\sigma_{3}\delta &=& \sigma_{1}\sigma_{2}.
\end{array} 
\end{equation}
Note that the action of $\delta$ on the five involutory generators 
\[\sigma_{1}\sigma_{2},\, \sigma_{1}\sigma_{2}\sigma_{3}\sigma_{1}^{-1},\, 
\sigma_{3}^{-1}\sigma_{1}\sigma_{2}\sigma_{3},\,  \sigma_{2}\sigma_{3},\,
\sigma_{1}\sigma_{2}\sigma_{3} \]
of $\Gamma({\cal P})$ can be illustrated by a Coxeter-type diagram on five nodes such that
$\delta$ corresponds to a symmetry of the diagram (see Figure~\ref{figthree} below for an example of a diagram on three nodes). In particular, if the nodes are given by the vertices and the center of a square diagram, with the central node representing $\sigma_{1}\sigma_{2}\sigma_{3}$
and the vertices representing the other four generators (in the order in which they are listed), then the action of $\delta$ on these generators corresponds to the $4$-fold rotational symmetry of the diagram. 
\smallskip

We now consider the mixing operation
\begin{equation}
\label{twistchirpc}
(\sigma_1\sigma_2, \sigma_1\sigma_2\sigma_3, \sigma_2\sigma_3, \delta) \mapsto 
(\delta, \sigma_1\sigma_2\delta^{-1}) =: (\kappa_{1},\kappa_{2})
\end{equation}
on $\overline{\Gamma }({\cal P})$, bearing in mind that the first three involutions on the left
side generate $\Gamma({\cal P})$. Set $\Lambda :=\langle \kappa_1,\kappa_2\rangle$. Then
$\Lambda = \overline{\Gamma }({\cal P})$; in fact, we have
\begin{equation}
\label{sigs}
\begin{array}{l}
\sigma_1 = \delta^2 \sigma_2 \sigma_3 = \delta^{-2} \sigma_2 \sigma_3 \delta \delta^{-1} = 
\delta^{-1} \sigma_1 \sigma_2 \delta^{-1} = \kappa_1^{-1} \kappa_2,\\
\sigma_2 = \sigma_2 (\sigma_3 \sigma_1 \sigma_2)^2 =  \delta \delta^{-1} \sigma_2 \sigma_3
\delta^2 \sigma_1 \sigma_2 =  (\delta \sigma_1 \sigma_2)^2 = 
(\sigma_1 \sigma_2 \delta^{-1})^{-2} = \kappa_2^{-2}, \\
\sigma_3 = \sigma_1 \sigma_2 \delta^2 = \sigma_1 \sigma_2 \delta^{-2} =
\kappa_2\kappa_1^{-1}.
\end{array}
\end{equation} 
Clearly, $\kappa_1$ is of order $4$, $\kappa_1 \kappa_2 = \delta \sigma_1 \sigma_2 \delta^{-1}$ is 
an involution, and $\kappa_1 \kappa_2^{-1} = \sigma_3^{-1}$ is of order $p$. Furthermore, 
\begin{equation}
\label{kapsquare}
\kappa_2^2 = \sigma_1 \sigma_2 \delta^{-1} \sigma_1 \sigma_2 \delta \delta^2 = 
\sigma_1 \sigma_2 \sigma_1 \sigma_2 \sigma_3 \sigma_1^{-1} \delta^2 = 
\sigma_3 \sigma_1^{-1} \sigma_1 \sigma_2 \sigma_3 = \sigma_2^{-1}, 
\end{equation}
and hence $\kappa_2^{2q} = \epsilon$. Moreover, $\kappa_2$ has order $2q$; in fact, $\kappa_2$,
being a duality, has even order, and $\sigma_2^{-1} = \kappa_2^2$ has order $q$. Hence, $\kappa_1,\kappa_2$ satisfy (at least) the following relations
\begin{equation}
\label{relch}
\kappa_1^4 = \kappa_2^{2q} = (\kappa_1 \kappa_2)^2 = (\kappa_1 \kappa_2^{-1})^p = \epsilon.
\end{equation}
These observations suggest that $\langle \kappa_1, \kappa_2 \rangle $ could be the
automorphism group of a chiral map. Indeed, this is true. 

\begin{biggy}
\label{impropthm}
Let $\cal P$ be an improperly self-dual chiral $4$-polytope of type $\{p,q,p\}$, let
$\overline{\Gamma }({\cal P})$ be its extended group, and let 
$\Lambda :=\langle \kappa_1,\kappa_2\rangle$, with $\kappa_1, \kappa_2$ as
in (\ref{twistchirpc}). Then $\Lambda$ is the automorphism group of a chiral $3$-polytope $\cal
M$ of type $\{4, 2q \hole p\}$, and $\Lambda = \overline{\Gamma }({\cal P})$.
\end{biggy}

\noindent {\bf Proof.}  We need to verify the intersection condition 
$\langle \kappa_1 \rangle \cap \langle \kappa_2 \rangle = \{\epsilon\}$ for $\Lambda$. 
Suppose that $\gamma \in \langle\kappa_1\rangle \cap \langle\kappa_2\rangle$. Then $\gamma
= \kappa_1^{j} = \delta^{j}$ for some $j \in \{0, 1, 2, 3\}$, and $\gamma = \kappa_{2}^{i}
= (\sigma_{1}\sigma_{2}\delta^{-1})^{i}$ for some $i$ with $0 \leq i < 2q$. It suffices to rule
out the possibility that $j=2$; then, by replacing $\gamma$ by $\gamma^{2}$, we have also ruled
out the cases $j=1,3$, and hence must have $j=0$.

Now suppose that $j=2$. Then $\delta^2 = \gamma = \kappa_{2}^i$ is an automorphism of $\cal P$,
and hence $i$ is even, $i = 2k$ (say). But then (\ref{kapsquare}) implies that
\[ \sigma_1\sigma_2\sigma_3 = \delta^2 = (\kappa_{2}^2)^{k} = \sigma_2^{-k}, \]
and hence $\sigma_{1} \in \langle\sigma_{2},\sigma_{3}\rangle$, a contradiction to the
intersection condition of $\Gamma({\cal P})$. Hence we cannot have $j=2$.

We now know that $\Lambda$ determines a chiral or regular map $\cal M$. It remains to show that
$\cal M$ cannot be regular. Suppose that $\cal M$ is regular. Then $\Lambda = 
\Gamma^{+}({\cal M})$. Let $\rho$ denote the $0^{\rm{th}}$ distinguished involutory generator of
$\Gamma({\cal M})$. Then conjugation in $\Gamma({\cal M})$ by $\rho$ gives
\[ \rho\kappa_1\rho = \kappa_1^{-1}, \;\; \rho\kappa_2\rho = \kappa_1^2 \kappa_2 .  \] 
Hence $\rho\delta\rho = \delta^{-1}$, and by (\ref{moredel}) and (\ref{sigs}) we have
\[ \rho\sigma_1\rho = \rho\kappa_1^{-1}\kappa_2\rho = \kappa_1 \kappa_1^2 \kappa_2 = 
\kappa_1^3 \kappa_2 = \kappa_1^{-1} \kappa_2 = \sigma_1 , \]
\[ \begin{array}{rl}
\rho\sigma_2\rho \,=&\!\!\! \rho\kappa_2^{-2}\rho = (\kappa_2^{-1}\kappa_1^{-2})^{2} =
\kappa_1 \kappa_2^{2} \kappa_1^{-1} = \delta (\sigma_1 \sigma_2 \delta^{-1})^2 \delta^{-1} =
\delta \sigma_1 \sigma_2 \delta^{-1} \sigma_1 \sigma_2 \delta^2 \\ 
=&\!\!\!\delta^2 (\delta^{-1} \sigma_1 \sigma_2 \delta) \delta^{2}\sigma_1\sigma_2\delta^{2}=
\delta^2 (\sigma_1 \sigma_2 \sigma_3 \sigma_1^{-1}) \sigma_1 \sigma_2 \sigma_3^2 = 
\sigma_2 \sigma_3^2, 
\end{array} \]
\[ \rho\sigma_3\rho = \rho\kappa_2 \kappa_1^{-1}\rho = \kappa_1^2 \kappa_2 \kappa_1 = 
\delta^2 \sigma_1 \sigma_2 \delta^{-1} \delta = \sigma_3^{-1} . \]
Thus $\rho$ determines an involutory group automorphism of $\Gamma({\cal P})$ that acts on
the generators $\sigma_{1},\sigma_{2},\sigma_{3}$ as indicated. Now it follows from \cite[Theorem 1]{chiral} that $\cal P$ must be regular, contrary to our assumptions; in fact, $\rho$ corresponds to
conjugation by the $3^{\rm{rd}}$ distinguished involutory generator of $\Gamma({\cal P})$.  Hence $\cal M$ must be chiral. This completes the proof.
\hfill $\Box$
\medskip

It is interesting to investigate the above construction in the context of {\em regular\/}
polytopes. Let $\cal P$ be a self-dual regular $4$-polytope, let
$\rho_0,\rho_1,\rho_2,\rho_{3}$ be the distinguished generators of $\Gamma({\cal P})$ with
respect to some base flag $\Phi$, and let $\sigma_{1},\sigma_{2},\sigma_{3}$ be the
corresponding distinguished generators of the rotation subgroup $\Gamma^{+}({\cal P})$, that
is, $\sigma_i := \rho_{i-1}\rho_i$ for $i=1,2,3$. As before, let $\delta$ be the duality (of
period $4$) defined by $\Phi\delta = \Phi^{0}$, and let $\kappa_{1},\kappa_{2}$ be as in
(\ref{twistchirpc}). Then the equations (\ref{delacts}), (\ref{moredel}) and (\ref{sigs})
carry over. Clearly, $\overline{\Gamma }({\cal P}) =
\langle\rho_0,\rho_1,\rho_2,\rho_3,\delta\rangle$. Since $\cal P$ is regular, we also have a
polarity $\omega$ that fixes $\Phi$. Then $\delta = \omega\rho_0\;(= \rho_{3}\omega)$. Once again, bear in mind here that dualities map a $j$-face in a flag to the $(3-j)$-face in the image flag.
Moreover, 
\[ \delta^{-1}\rho_{j}\delta = \rho_{0}\omega\rho_{j}\omega\rho_{0} = 
\rho_{0}\rho_{3-j}\rho_{0} ,\]
so $\delta$ acts on the generators of $\Gamma({\cal P})$ as follows:
\[ \delta^{-1}\rho_0\delta = \rho_3, \;\; 
\delta^{-1}\rho_1\delta = \rho_2, \;\;
\delta^{-1}\rho_2\delta = \rho_0\rho_1\rho_0,\;\;
\delta^{-1}\rho_3\delta = \rho_0.  \]

Define $\rho:=\rho_{3}$. (Recall that the supposed group automorphism of $\Gamma({\cal P})$ in the proof of Theorem~\ref{impropthm} was determined by conjugation with $\rho_{3}$.)  Then, using $\delta = \omega\rho_0$, the generators of $\Gamma({\cal P})$ can be expressed in terms of $\kappa_1,\kappa_2,\rho$ as follows:
\[ \rho_0 = \rho \kappa_1^{2}\, (= \rho \delta^2),\;\; 
\rho_1 =  \rho\kappa_1\kappa_2,\;\;
\rho_2 = \rho \kappa_1 \kappa_2^{-1},\;\;
\rho_3 = \rho. \]

Now the mixing operation (\ref{twistchirpc}) can be rewritten in terms of the generators $\rho_i$ as follows:
\begin{equation}
\label{twistregpc}
(\rho_0, \rho_1, \rho_2, \rho_3, \delta) \mapsto ( \rho_3, \rho_3 \delta, \rho_1 ) 
= (\rho, \rho \kappa_1 , \rho \kappa_1 \kappa_2 ) .
\end{equation}
(Recall here that, in the proof of Theorem~\ref{impropthm}, the element $\rho$ acted on $\cal M$ like the $0^{\rm{th}}$ distinguished generator.)  Note that $(\rho\kappa_1)^{2} = \epsilon$, so the generators of the corresponding rotation subgroup are indeed $\kappa_1,\kappa_2$. Moreover, $\rho_{3}\delta = \omega$, so the operation (\ref{twistregpc}) is simply the dual version of (\ref{twistpcp}) (obtained by conjugating the new generators with $\omega$). 

We summarize our observations in 

\begin{riggy}
\label{rem}
Let $\cal P$ be a self-dual regular $4$-polytope of type $\{p,q,p\}$, let
$\overline{\Gamma }({\cal P})$ be its extended group, and let 
$\Lambda :=\langle \kappa_1,\kappa_2\rangle$, with $\kappa_1, \kappa_2$ as in (\ref{twistchirpc}). Then $\Lambda$ is the rotation subgroup of a regular $3$-polytope $\cal M$ of type $\{4, 2q \hole p\}$, whose full automorphism group is $\overline{\Gamma }({\cal P})$. In particular, $\cal M$ is the 
regular map described in \rm{\cite{peter-egon}} and defined by (\ref{twistpcp}).
\end{riggy}

We now give some examples of chiral maps that are associated with improperly self-dual locally toroidal chiral 4-polytopes.  For notation we refer to \cite[\S 2]{chirality}.
\medskip

\begin{figure}[hbt]
\begin{center}
\begin{picture}(150,160)
\includegraphics[width=2in]{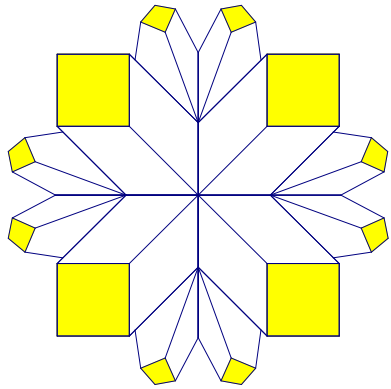}
\end{picture}
\caption{\em A chiral map of type $\{4,8 \hole 4\}$}
\label{figone}
\end{center}
\end{figure}

\noindent {\bf Example 1}\hskip.1in  The universal chiral $4$-polytope ${\cal P} = 
\{ \{4,4\}_{(1,3)},\{4,4\}_{(1,3)}\}$ is improperly self-dual and has an automorphism group of order $2000$ (see \cite[p.132]{self-duality}). The induced chiral map $\cal M$ of Theorem~\ref{impropthm}, with an automorphism group of order $4000$, is of type  $\{4,8 \hole 4\}$.  Eight quadrangular faces meet at each vertex of the map. The local structure is seen in Figure~\ref{figone}, where the $4$-gonal holes of the map are highlighted (they are not faces of the map). In a sense, the set of quadrangular faces splits into blocks of $4$, each associated with a hole and corresponding to the block of rectangular faces (in the mantle) of a $4$-gonal prism; in the figure, the central vertex is common to the four prisms that are associated with the four highlighted square-shaped holes. It should be noted that the valency of each vertex is $8$, but some edges in the figure are hidden. In addition we observe that the polytope $\cal P$ has no polarity and hence the group of the map $\cal M$ cannot be generated by involutions. 
\medskip

\medskip

\noindent {\bf Example 2}\hskip.1in  The universal chiral $4$-polytope ${\cal P} = 
\{ \{6,3\}_{(1,2)},\{3,6\}_{(2,1)}\}$ is improperly self-dual and has an automorphism group of order $20160$. The polytope has (at least) two quotients, with groups of orders $10080$ and $5040$ (in the latter case, isomorphic to the symmetric group $S_7$), which are also improperly self-dual chiral polytopes of type $\{6,3,6\}$ with toroidal facets $\{6, 3\}_{(1,2)}$ and vertex-figures $\{3,6\}_{(2,1)}$;  both can be derived from $\cal P$ by identifying vertices separated by a certain number of steps along Petrie polygons. The polytope $\cal P$ and its two quotients have Petrie polygons of length $28$, $14$ and $7$, respectively.  Now all three polytopes possess a polarity (see \cite[Theorem 3.4]{self-duality}), so their extended groups are generated by involutions. The automorphism groups of the three induced maps $\cal M$, of type $\{4,6 \hole 6\}$, are thus generated by involutions.  The local structure of the maps, with hexagonal holes highlighted, can be seen in Figure~2; now the set of quadrangular faces splits into blocks of $6$, each associated with a hole and corresponding to a $6$-gonal prism. 
\medskip

\begin{figure}[hbt]
\begin{center}
\begin{picture}(150,160)
\includegraphics[width=2in]{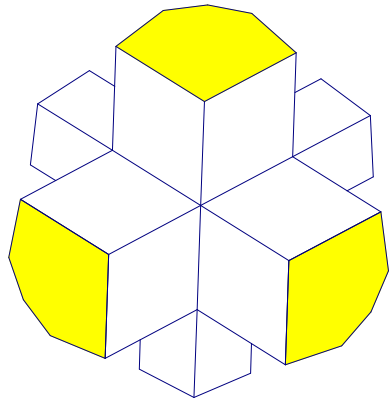}
\end{picture}
\caption{\em A chiral map of type $\{4,6 \hole 6\}$}
\label{figtwo}
\end{center}
\end{figure}

\medskip

It should be noted that, for improperly self-dual chiral $4$-polytopes, the elements $\pi_L = \sigma_1 \sigma_3$ and $\pi_R = \sigma_1 \sigma_3^{-1} = \delta^{-1} \sigma_3 \sigma_1 \delta$ are conjugate in $\overline{\Gamma }({\cal P})$ (since $\sigma_3\sigma_1$ and $\sigma_1\sigma_3$ are conjugate in $\Gamma({\cal P})$) and the left- and right-handed Petrie polygons are thus of the same length. 

\section{Regular maps from properly self-dual chiral polytopes}
\label{chprop}

As before, let $\cal P$ be a chiral $4$-polytope of type $\{p,q,p\}$, and let $\sigma_1,\sigma_2,\sigma_3$ be the distinguished generators of $\Gamma (\cal P)$ defined with respect to a base flag $\Phi$. We now assume that $\cal P$ is properly self-dual.  Then $\cal P$ admits a (unique) polarity $\omega$ which fixes $\Phi$, and
\begin{equation}
\label{omsig}
\omega \sigma_1 \omega = \sigma_3^{-1} ,\quad
\omega \sigma_2 \omega = \sigma_2^{-1} 
\end{equation}
in $\overline{\Gamma }({\cal P})$ (see \cite[Theorem 3.1]{self-duality}). It follows that
\begin{equation}
\omega \sigma_1 \sigma_2 \omega = \sigma_2 \sigma_3, \quad 
\omega \sigma_1 \sigma_2 \sigma_3 \omega = \sigma_1 \sigma_2 \sigma_3 ,
\end{equation}
so that (under conjugation) $\omega$ fixes $\sigma_1 \sigma_2 \sigma_3$ and interchanges $\sigma_1 \sigma_2$ with $\sigma_2 \sigma_3$. Hence $\omega$ induces the twisting operation 

\begin{equation}
\label{twistprop}
(\sigma_1 \sigma_2 \sigma_3, \sigma_1 \sigma_2, \sigma_2 \sigma_3, \omega) 
\mapsto (\sigma_1 \sigma_2 \sigma_3, \sigma_1 \sigma_2, \omega) =:
(\tau_0,\tau_1,\tau_2)
\end{equation}

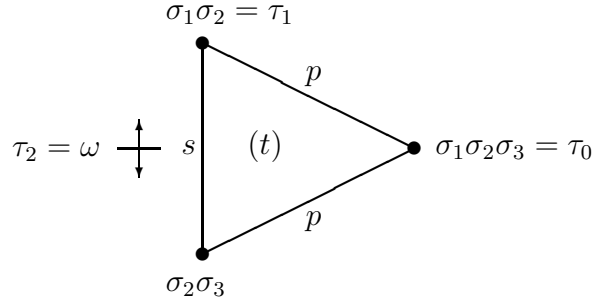
\begin{figure}[htb]
\centering
\begin{picture}(100,130)
\thicklines
\multiput(20,30)(0,80){2}{\circle*{5}}
\put(100,70){\circle*{5}}
\put(20,30){\line(2,1){80}}
\put(20,110){\line(2,-1){80}}
\put(20,30){\line(0,1){80}}
\put(-52,67){$\tau_{2}=\omega$}
\thinlines
\put(-4,70){\vector(0,1){11}}
\put(-4,70){\vector(0,-1){11}}
\put(-12,70){\line(1,0){16}}
\put(6,16){$\sigma_2\sigma_3$}
\put(6,119){$\sigma_1\sigma_2 = \tau_1$}
\put(108,68){$\sigma_1\sigma_2\sigma_3 = \tau_0$}
\put(12,68){$s$}
\put(59,96){$p$}
\put(59,41){$p$}
\put(37,68){$(t)$}
\end{picture}
\caption{\em The polarity $\omega$ acting on $\Gamma({\cal P})$.}
\label{figthree}
\end{figure}

\noindent on $\Gamma({\cal P})$ shown in Figure~\ref{figthree}; the notation used in the figure is explained in more detail below. Note that each $\tau_i$ is an involution, and that $\tau_0 $ and $\tau_2 $ commute. Furthermore,
\[ \begin{array}{l}
\sigma_1 = \sigma_1 (\sigma_2 \sigma_3)^2 = 
\sigma_1 \sigma_2 \sigma_3 \omega \sigma_1 \sigma_2 \omega = 
\tau_0 \tau_2 \tau_1 \tau_2 = \tau_2 \tau_0 \tau_1 \tau_2 ,\\
\sigma_2 = \sigma_1^{-1} \tau_1 = \tau_2 \tau_1 \tau_2 \tau_0 \tau_1 ,\\
\sigma_3 = \sigma_1 \sigma_2 \tau_0 = \tau_1 \tau_0 .
\end{array} \]
Hence $\langle \tau_0,\tau_1,\tau_2 \rangle = \overline{\Gamma }({\cal P})$, and 
$\tau_{0}\tau_{1} = \sigma_{3}^{-1}$ is of order $p$. 

\begin{biggy}
\label{propthm}
Let $\cal P$ be a properly self-dual chiral $4$-polytope of type $\{p,q,p\}$, let
$\overline{\Gamma }({\cal P})$ be its extended group, and let 
$\Lambda :=\langle \tau_0,\tau_1,\tau_2 \rangle$, with $\tau_0,\tau_1,\tau_2$ as
in (\ref{twistprop}). Then $\Lambda$ is the automorphism group of a regular $3$-polytope 
$\cal M$ of type $\{p,2s\}$, where $s$ is the length of the left-Petrie polygons of $\cal P$.
Moreover, $\Lambda = \overline{\Gamma }({\cal P})$.
\end{biggy}

\noindent {\bf Proof.}  First observe that
\[ (\tau_1\tau_2)^{2} = (\sigma_1\sigma_2 \omega)^2 = \sigma_1 \sigma_2^2 \sigma_3 =
\sigma_2^{-1} \sigma_1^{-1} \sigma_3^{-1} \sigma_2^{-1} 
= \sigma_2^{-1} \sigma_2 \sigma_3 \sigma_1  = \sigma_3 \sigma_1 , \] 
which is conjugate to $\pi_{L} := \sigma_1\sigma_3$ and hence is of order $s$. Then, since the order of a duality is even, $\tau_1\tau_2$ must have order $2s$. 

It remains to verify the intersection condition for $\Lambda$.  Suppose that
$\gamma \in \langle\tau_0,\tau_1\rangle  \cap \langle\tau_1,\tau_2\rangle$.  By multiplying $\gamma$ by $\tau_1$ if need be, we may assume that $\gamma = (\tau_1\tau_0)^{j} = \sigma_3^j$ for some $j$, so that $\gamma$ fixes the base vertex $F_0$ of $\Phi$. We claim that $\gamma=\epsilon$. Now, since
$\gamma \in \langle \tau_1,\tau_2\rangle 
= \langle \sigma_1 \sigma_2,\sigma_2\sigma_3\rangle  \ltimes  \langle\omega\rangle$ 
and $\gamma \in \Gamma (\cal P)$, we also have
\[ \gamma = (\sigma_1 \sigma_2^2 \sigma_3)^k (\sigma_2 \sigma_3)^l 
= (\sigma_3 \sigma_1)^k (\sigma_2 \sigma_3)^l \] 
for some integer $k$ and $l= 0,1$. If $l=0$, then $\gamma$ must shift the base vertex $F_0$ $k$ steps along a left-Petrie polygon (the image under $\sigma_3^{-1}$ of the left-Petrie polygon associated with $\pi_L$), so that necessarily $k=0$; it follows that $\gamma = \epsilon$, as required.  On the other hand, the case $l=1$ cannot occur. In fact, when $l=1$, the element $(\sigma_3 \sigma_1)^k = \gamma (\sigma_2 \sigma_3)^{-1}$ also fixes $F_0$, so that again $k=0$ and now $\gamma = \sigma_2 \sigma_3$;  but from $\gamma = \sigma_3^j$ we then obtain $\sigma_2 = \sigma_3^{j-1}$, a contradiction.  

It follows that $\Lambda$ is the group of a regular $3$-polytope, which then must be of type $\{p,2s\}$.
\hfill $\Box$
\medskip\smallskip

We can also determine the lengths of the Petrie polygons and $2$-zigzags of $\cal M$.

\begin{riggy}
\label{remtwo}
Let $\cal P$ and $\cal M$ be as in Theorem~\ref{propthm}, and let $t$ be the 
length of the right-Petrie polygons of $\cal P$. Then $\cal M$ is of type
$\{p,2s\}_{2t,q}$; that is, the length of the Petrie polygons of $\cal M$ is $2t$,
and the length of the $2$-zigzags of $\cal M$ is $q$. 
\end{riggy}

The proof is straightforward.  In fact,
\[ (\tau_0 \tau_1 \tau_2 )^2 = (\sigma_3^{-1} \omega )^2 = \sigma_3^{-1} \sigma_1 \]
is conjugate to $\pi_R = \sigma_1 \sigma_3^{-1}$ and hence has order $t$.  On the 
other hand, since dualities have even order, the duality $\tau_0\tau_1\tau_2$ then
must have order $2t$. Finally, $\tau_0 (\tau_1 \tau_2)^2$ is conjugate to~$\sigma_2$. 
\bigskip

The diagram shown in Figure~\ref{figthree} exhibits the group $\Gamma({\cal P})$ as a quotient of the abstract group $\Gamma^{3}(s,t;p,p)$ described in \cite[Sect. 9D]{arp}.  As usual, the mark on a branch of the diagram is the period of the product of the two generators at its ends. The interior mark, $t$, represents the extra relation 
\[ (\sigma_1\sigma_2 \cdot \sigma_1\sigma_2\sigma_3 \cdot \sigma_2\sigma_3 \cdot
\sigma_1\sigma_2\sigma_3)^{t} = \epsilon , \]
which holds in $\Gamma({\cal P})$ because of
\[ \sigma_1\sigma_2 \cdot \sigma_1\sigma_2\sigma_3 \cdot \sigma_2\sigma_3 \cdot
\sigma_1\sigma_2\sigma_3 
\;=\; (\sigma_1\sigma_2)^{2}\sigma_3\sigma_{1}^{-1} (\sigma_1\sigma_2\sigma_3)^{2} 
\;=\; \sigma_3\sigma_{1}^{-1} \;=\; \pi_R^{-1} .\] 
Note that the mark $t$ is placed in parentheses to indicate the fact that the other analogous products of four generators, with one repeated, will not in general have the same order $t$ (see \cite[p.320]{arp} for more details).
\bigskip

\noindent {\bf Example 3}\  The universal chiral $4$-polytope ${\cal P} = 
\{\{3,6\}_{(1,2)},\{6,3\}_{(1,2)}\}$ is properly self-dual. Its automorphism group, of order $672$, is the group $\rm{L}_{2}^{\langle \pm 1\rangle}(\mathbb{Z}_7)$ of all $2\times 2$ matrices of determinant $\pm 1$ over $\mathbb{Z}_7$ (see \cite{chirality}). The induced regular map $\cal M$ is of type $\{3,16\}$ (more exactly, of type $\{3,16\}_{28,6}$) and has $42$ vertices and $224$ faces; its automorphism group is 
${L}_{2}^{\langle \pm 1\rangle}(\mathbb{Z}_7) \ltimes C_2$ of order $1344$. The quotient ${\cal P}'$ of $\cal P$ by the central subgroup $\langle \pm I\rangle$ (with $I$ the identity matrix) of  $\rm{L}_{2}^{\langle \pm 1\rangle}(\mathbb{Z}_7)$ again has facets $\{3,6\}_{(1,2)}$ and vertex-figures $\{6,3\}_{(1,2)}$, and yields a regular map of type $\{3,8\}$ with $42$ vertices, $112$ faces, and automorphism group $\rm{PGL_{2}(7)} \ltimes C_2$ of order $672$. The latter has Petrie polygons of length $14$, and doubly covers $\{3,8\}_7$, the dual of the Petrial of the Klein map $\{7,3\}_8$ (see \cite{cm}).

\section{Involutory generators for chiral maps}

In the construction of polytopes from those of higher rank we often encounter the
question whether or not the automorphism group or a certain subgroup is generated by
involutions. We know that the group of a regular polytope of any rank $n$ is a
C-group and comes with a  distinguished set of generators consisting of involutions.
The group of a chiral polytope $\cal P$ of rank $n\geq 4$ is also generated by
involutions (see \cite[p. 496]{chiral}), although its distinguished generators $\sigma_1,\ldots,\sigma_{n-1}$ are not involutory; for example, when $n=4$, the elements
$\sigma_1\sigma_2\sigma_3,\sigma_1\sigma_2,\sigma_2\sigma_3$ are involutory generators. However, this is no longer true for rank $3$, as shown by the chiral map described in Example~1
(as well as the toroidal maps below).  On the other hand, many chiral $3$-polytopes (see Example~2) do have the desired property.  In this Section we briefly discuss obstructions for the existence of involutory generators for rank~$3$.

Let $\cal P$ be a chiral $3$-polytope of type $\{p,q\}$ with group 
$\Gamma({\cal P}) = \langle\sigma_{1},\sigma_{2}\rangle$. Let
$\tau:=\sigma_1\sigma_2$, and let
\begin{equation}
\label{ntau}
N({\tau}) := \langle \varphi^{-1}\tau\varphi \mid 
\varphi\in\Gamma({\cal P})\rangle
\end{equation} 
(that is, $N({\tau})$ is the normal closure of $\tau$ in $\Gamma({\cal P})$). Then
$\tau$ is the ``half-turn" about the base edge (in $\Phi$), and since $\Gamma({\cal
P})$ is edge-transitive, $N({\tau})$ is the group generated by all the half-turns
about  edges of $\cal P$. Clearly, since $\tau$ is an involution, $N({\tau})$ is
generated by involutions (although their number may not be finite if $\cal P$ is
infinite). Now, since $\Gamma({\cal P}) = \langle\tau,\sigma_{2}\rangle =
\langle\sigma_{1},\tau\rangle$, we have 
\[ \Gamma({\cal P}) = N({\tau}) \cdot \langle\sigma_{2}\rangle =
N({\tau}) \cdot \langle\sigma_{1}\rangle \]
(a product of subgroups), where the second factors are cyclic of orders $q$ or
$p$, respectively. In particular, $\Gamma({\cal P})/N({\tau})$ is a cyclic group,
whose order divides $p$ and $q$ if the latter are finite. 

\begin{liggy}
\label{lemmaone}
Let $\Gamma$ be any group, and let $N$ be a normal subgroup of $\Gamma$ such that 
$\Gamma/N$ is cyclic. If $\Gamma$ is generated by involutions, then $N$ has index $1$ or $2$ in $\Gamma$.
\end{liggy}

\noindent {\bf Proof.}  Suppose that $N \neq \Gamma$. If $\Gamma$ is generated by
involutions, then so is its quotient $\Gamma/N$. A non-trivial cyclic group
contains an involution only if it has finite even order, and it is generated by
involutions only if its order is $2$. This proves the lemma.
\hfill $\Box$
\medskip

The following proposition follows immediately from Lemma~\ref{lemmaone}, applied with
$\Gamma=\Gamma({\cal P})$ and $N=N({\tau})$.

\begin{piggy}
\label{propone}
Let $\cal P$ be a chiral $3$-polytope of type $\{p,q\}$ with group 
$\Gamma({\cal P}) = \langle\sigma_{1},\sigma_{2}\rangle$, and let $N({\tau})$ be as
in (\ref{ntau}). \\
{\rm(a)} If $\Gamma({\cal P})$ is generated by involutions, then $N({\tau})$ has
index $1$ or $2$ in $\Gamma({\cal P})$.\\
{\rm(b)} If $\Gamma({\cal P})=N({\tau})$ or $N({\tau}) \ltimes C_2$, then
$\Gamma({\cal P})$ is generated by involutions.  
\end{piggy}

The proposition basically says that, if $\cal P$ is a chiral $3$-polytope with
involutory generators, then in fact the half-turns about edges provide a
``nice" generating set of involutions for a subgroup of very small index, at most
$2$. The  second part of Proposition~\ref{propone} is a partial converse of the
first. Note that, if $N({\tau})$ has index $2$ and $p\equiv 2\, (\bmod\;4)$ or
$q\equiv 2\,(\bmod\;4)$, respectively, then $\Gamma({\cal P}) =
N({\tau})\ltimes\langle\sigma_{1}^{\,p/2}\rangle$ or
$\Gamma({\cal P}) = N({\tau})\ltimes \langle\sigma_{2}^{\,q/2}\rangle$; in any case,
$\Gamma({\cal P}) = N({\tau}) \ltimes C_2$ and hence
$\Gamma({\cal P})$ is generated by involutions.

We remark that there is an exact analogue to Proposition~\ref{propone} for the
rotation subgroup $\Gamma^{+}({\cal P})$ of a regular $3$-polytope. The rotation
subgroups of regular polytopes of rank at least $4$ are always generated by
involutions.

We conclude this Section with a brief discussion of the chiral maps on the torus.  None of these maps $\{4,4\}_{(b,c)}$, $\{3,6\}_{(b,c)}$ and $\{6,3\}_{(b,c)}$ has a group generated by involutions. 

Consider, for example, ${\cal P} := \{4,4\}_{(b,c)}$. Then $\cal P$ is a quotient of the regular plane tessellation $\{4,4\}$ with vertex-set $\mathbb{Z}^2$. Let $[4,4]^{+}$ denote the rotation subgroup of $\{4,4\}$, with distinguished generators $\sigma_1,\sigma_2$, let $T$ be its translation subgroup generated by the translations $\tau_{1},\tau_{2}$ by the vectors $(1,0)$ and $(0,1)$, respectively, and let $T_{(b,c)}$ denote the subgroup of $T$ defined by
\[ T_{(b,c)} = \langle \tau_{1}^{b}\tau_{2}^{c},\tau_{1}^{-c}\tau_{2}^{b} \rangle . \]
Then $[4,4]^{+} = T \ltimes \langle \sigma_2 \rangle \cong \mathbb{Z}^2 \ltimes C_4$ and
\[ [4,4]_{(b,c)} := \Gamma(\{4,4\}_{(b,c)}) \cong [4,4]^{+} / T_{(b,c)} 
\cong (T/T_{(b,c)}) \ltimes C_4 \]
(see \cite[\S 8.3]{cm}); note here that $\langle \sigma_2 \rangle \cap T_{(b,c)}$ is trivial, so the semi-direct product structure is preserved under taking the quotient. But now we can simply appeal to Lemma~\ref{lemmaone}, applied with $\Gamma = [4,4]_{(b,c)}$ and $N = T/T_{(b,c)}$. It follows that
$[4,4]_{(b,c)}$ is not generated by involutions, since the index of $N$ is $4$.

Similar arguments also apply to the chiral maps $\{3,6\}_{(b,c)}$. Now we must take 
\[ T_{(b,c)} := \langle \tau_{1}^{b+c}\tau_{2}^{c},\tau_{1}^{-c}\tau_{2}^{b} \rangle  \]
(see \cite[\S 8.4]{cm}).  The dual map $\{6,3\}_{(b,c)}$ has the same group as $\{3,6\}_{(b,c)}$, and hence is also not generated by involutions. 

Moreover, the same arguments establish that the rotation subgroups for the regular torus maps 
$\{4,4\}_{(b,c)}$, $\{3,6\}_{(b,c)}$ and $\{6,3\}_{(b,c)}$ (with $c=0$ or $b=c$) cannot be generated by involutions.

\section{Concluding remarks}
\label{conrem}

In his 1948 paper \cite{config}, Coxeter enumerated the chiral maps on the torus
and showed that they fall into the three infinite families of maps of type $\{4,4\}$,
$\{3,6\}$ or $\{6,3\}$.  Garbe~\cite{garbe} proved that there are no chiral maps on
surfaces of genus $2, 3, 4, 5$ or $6$.  While there are many regular maps known to
exist, the occurrence of chiral maps is rather sporadic. Several attempts have been
made to construct chiral maps, notably by Sherk~\cite{sherk} and Wilson~\cite{doro-wilson, wilson}.  More recently, in \cite{conder}, Conder and Dobcs\'{a}nyi have enumerated all regular and chiral maps on surfaces up to genus $15$. However, of the sixteen chiral maps listed, only thirteen are in fact abstract polyhedra, that is, their automorphism groups satisfy the intersection condition. (Every abstract polyhedron yields a map on a surface, but not vice versa. A map may be ``too small" to be polytopal.  An example is the torus map $\{4,4\}_{(1,0)}$, which has a single vertex, two edges and a single face.) The construction described in this paper can be used to provide new examples of chiral maps. They all are related to improperly self-dual chiral $4$-polytopes. Unfortunately, at present, while many examples (including infinite families) of properly self-dual chiral $4$-polytopes are known (see, for example, \cite{chirality}), not many examples of improperly self-dual chiral $4$-polytopes can be found in the literature. 

\medskip
We are grateful to Barry Monson for his help in verifying the existence of certain polytopes using computer enumerations. We also thank the anonymous referee for a number of helpful comments.

\vskip.2in
\noindent
{\em Addresses\/}:\\[.1in]
{\em Isabel Hubard, York University, Toronto, ON, Canada M3J 1P3.}\\[.03in] 
{\em Egon Schulte, Northeastern  University, Boston, MA 02169, USA.}\\[.03in] 
{\em Asia Ivi\'{c} Weiss, York University, Toronto, ON, Canada M3J 1P3.}

\end{document}